\documentclass[12pt,english]{article}
\usepackage[latin9]{inputenc}
\usepackage{geometry}
\usepackage{babel}
\usepackage{amstext}
\usepackage{amssymb}
\usepackage{wasysym}
\usepackage[unicode=true,pdfusetitle,
 bookmarks=true,bookmarksnumbered=false,bookmarksopen=false,
 breaklinks=false,pdfborder={0 0 1},backref=false,colorlinks=false]
 {hyperref}

\makeatletter

\providecommand{\tabularnewline}{\\}

\date{}

\makeatother

\begin{document}

\title{An Achievement Game on a Cycle}

\author{Eero R\"{a}ty\thanks{Centre for Mathematical Sciences, Wilberforce Road, Cambridge CB3
0WB, UK, epjr2@cam.ac.uk}}
\maketitle
\begin{abstract}
Consider the following game played by Maker and Breaker on the vertices
of the cycle $C_{n}$, with first move given to Breaker. The aim of
Maker is to maximise the number of adjacent pairs of vertices that
are both claimed by her, and the aim of Breaker is to minimise this
number. The aim of this paper is to find this number exactly for all
$n$ when both players play optimally, answering a related question
of Dowden, Kang, Mikala\v{c}ki and Stojakovi\'{c}. 
\end{abstract}

\section{Introduction}

Consider the following game, called the 'Toucher-Isolator' game on
a graph $G$, introduced by Dowden, Kang, Mikala\v{c}ki and Stojakovi\'{c}
\cite{key-7}. The two players, Toucher and Isolator, claim edges
of $G$ alternately with Toucher having the first move. Let $t\left(G\right)$
be the number of vertices that are incident to at least one of the
edges claimed by Toucher. The aim of Toucher is to maximise $t\left(G\right)$
and the aim of Isolator is to minimise $t\left(G\right)$. Hence this
is a 'quantitative' Maker-Breaker type of game. 

For given $G$, let $u\left(G\right)$ be the number of isolated vertices,
i.e. the number of vertices that are not incident to any of the edges
claimed by Toucher at the end of the game when both players play optimally.
Dowden, Kang, Mikala\v{c}ki and Stojakovi\'{c} gave bounds \cite{key-7} for
the size of $u\left(G\right)$ for general graphs $G$, and studied
some particular examples as well which included cycles and paths.
In particular, they proved that

\[
\frac{3}{16}\left(n-3\right)\le u\left(C_{n}\right)\leq\frac{n}{4}
\]
and 

\[
\frac{3}{16}\left(n-2\right)\leq u\left(P_{n}\right)\leq\frac{n+1}{4},
\]
where $C_{n}$ is a cycle with $n$ vertices and $P_{n}$ is a path
with $n$ vertices. 

Note that these bounds imply that the asymptotic proportion of untouched
vertices is between $\frac{3}{16}$ and $\frac{1}{4}$ in both cases.
Dowden, Kang, Mikala\v{c}ki and Stojakovi\'{c} asked what the correct asymptotic
proportion of untouched vertices is, and suggested that the correct
answer could be $\frac{1}{5}$. In this paper we prove that this is
the correct asymptotic proportion, and in fact we give the exact values
of $u\left(C_{n}\right)$ and $u\left(P_{n}\right)$ for all $n$.
\\

\textbf{Theorem 1. }When $G=C_{n}$ and both players play optimally,
there will be $\left\lfloor \frac{n+1}{5}\right\rfloor $ untouched
vertices. \\

\textbf{Theorem 2. }When $G=P_{n}$ and both players play optimally,
there will be $\left\lfloor \frac{n+4}{5}\right\rfloor $ untouched
vertices. \\

Although this paper is self-contained, for general background on Maker-Breaker
type games, see Beck \cite{key-2}. There are many other papers dealing
with achievement games on graphs, see e.g.\ \cite{key-3,key-4,key-5}.

For convenience we work on a 'dual version' of these games. Consider
a game played on the vertices of a cycle $C_{n}$ with two players
Maker and Breaker claiming vertices in alternating turns, with first
move given to Breaker. For this game, define the \textit{score} to
be the number of adjacent pairs of vertices claimed by Maker on the
cycle. It is easy to see that this game is identical to the Toucher-Isolator
game played on $C_{n}$, with Maker corresponding to Isolator and
Breaker corresponding to Toucher. Indeed, this follows from the fact
that claiming adjacent pairs of vertices on the dual game corresponds
to claiming two edges whose endpoints meet in the original game, which
is precisely the same as isolating the vertex where they meet. 

When considering the dual version for the path, we have to be a bit
more careful due to irregular behaviour at the endpoints. For that
reason it turns out to be useful to define three different games which
essentially only differ at the endpoints of the path. Firstly define
a game $F\left(n\right)$ played on the elements of $\left\{ 1,\dots,n\right\} $
with two players Maker and Breaker claiming elements in alternating
turns with first move given to Maker. For this game, define the \textit{score}
to be the number of pairs $\left\{ i,i+1\right\} $ such that both
$i$ and $i+1$ are claimed by Maker, and as usual Maker is aiming
to maximise this score and Breaker is aiming to minimise this score.
Let $\alpha\left(n\right)$ be the score attained when both Maker
and Breaker play optimally.

Similarly as with $C_{n}$, the game $F\left(n-1\right)$ and the
Toucher-Isolator game on $P_{n}$ have a strong relation with Maker
corresponding to Isolator and Breaker corresponding to Toucher. However
unlike with the game on the cycle, they are not exactly the same game
as first and last vertex can be isolated by claiming only the first
or last edge respectively, and also because Toucher has first move
in the Toucher-Isolator game whereas Maker has first move in the game
$F\left(n-1\right)$. 

We also define the games $G\left(n\right)$ and $H\left(n\right)$
both played on $\left\{ 1,\dots,n\right\} $, with players Maker and
Breaker claiming elements in alternating turns with first move given
to Maker. On $G\left(n\right)$, we increase the score by one for
each pair $\left\{ i,i+1\right\} $ with both $i$ and $i+1$ claimed
by Maker, and the score is also increased by $1$ if Maker claims
the element $1$. In a sense, this can be viewed as a game on the
board $\left\{ 0,\dots,n\right\} $ with $0$ assigned to Maker initially.
Similarly on $H\left(n\right)$, we increase the score by one for
each pair $\left\{ i,i+1\right\} $ with both $i$ and $i+1$ claimed
by Maker, and additionally score is increased by $1$ for claiming
either of the elements $1$ or $n$. Again, this can be viewed as
a game on the board $\left\{ 0,\dots,n+1\right\} $ with both $0$
and $n+1$ assigned to Maker initially. Define $\beta\left(n\right)$
and $\gamma\left(n\right)$ be the scores of these games when both
players play optimally. 

The idea behind defining these games is the following. If $B$ is
a game of the form $F\left(n\right)$, $G\left(n\right)$ or $H\left(n\right)$,
and if Breaker plays her first move adjacent to Maker's first move,
then the board $B$ splits into two disjoint boards, which are of
the form $F\left(m\right)$, $G\left(m\right)$ or $H\left(m\right)$
- however, note that these two boards are not in general of the same
form, and not necessarily of the same form as the original board.
Hence it turns out to be useful to analyse all of these games at the
same time. 

Consider the dual game played on $C_{n}$ with vertex set $\left\{ 1,\dots,n\right\} $,
and recall that in this dual version the first move is given to Breaker.
By the symmetry of the cycle, we may assume that Breaker claims $n$
on her first move. Hence after this first move, the available winning
lines that can increase the score are $\left\{ 1,2\right\} ,\dots,\left\{ n-2,n-1\right\} $.
These are exactly the winning lines of the game $F\left(n-1\right)$,
and since Maker has the next move it follows that the subsequent game
is equivalent to the game $F\left(n-1\right)$. Hence $u\left(C_{n}\right)=\alpha\left(n-1\right)$,
and thus it suffices to find the value of $\alpha\left(n\right)$
for all $n$. 

In order to analyse the Toucher-Isolator game on $P_{n}$, define
the game $H_{b}\left(n\right)$ in exactly the same way as $H\left(n\right)$,
but with the first move given to Breaker, and let $\gamma_{b}\left(n\right)$
be the score of this game when both players play optimally. It is
easy to see that the Toucher-Isolator game and $H_{b}\left(n\right)$
are equivalent in the same sense as the Toucher-Isolator game on $C_{n}$
and $F\left(n-1\right)$ are. Hence it follows that $u\left(P_{n}\right)=\gamma_{b}\left(n\right)$,
and thus it suffices to find the value of $\gamma_{b}\left(n\right)$. 

We start by focusing on $F\left(n\right)$ and finding the value $\alpha\left(n\right)$.
Since Maker is trying to maximise the score, it seems sensible for
her to start by claiming some suitably chosen $i$, and then trying
to claim as long block of consecutive elements as possible. As long
as $i\not\in\left\{ 1,n\right\} $, she can certainly guarantee a
block of length at least 2. Now suppose she has claimed a block of
length $t$, and she cannot proceed in this way. This means that Breaker
must have claimed the points next to the endpoints of this block (or
one of the endpoints is $1$ or $n$). Removing this block, together
with the endpoints Breaker has claimed, leaves a path with $n-t-1$
elements containing at most $t-1$ elements claimed by Breaker, and
no elements claimed by Maker. 

This motivates the definition of the following game, which can be
viewed as a delayed version of $F\left(n\right)$. Let $F\left(n,k\right)$
be the game played on $\left\{ 1,\dots,n\right\} $, where at the
start of the game Breaker is allowed to claim $k$ points, and then
the players claim elements alternately, with the score defined in
the same way as for $F\left(n\right)$. Thus $F\left(n\right)$ and
$F\left(n,0\right)$ are identical games. 

Let $\alpha\left(n,k\right)$ be the score attained when both players
play optimally. It turns out that by following the strategy described
above with a suitable choice of the initial move, we can prove a good
enough lower bound for $\alpha\left(n,k\right)$, and almost same
argument also works for $\gamma_{b}\left(n\right)$. 

One can observe from the proof of the lower bound of $\alpha\left(n,k\right)$
that allowing Maker to have multiple 'long blocks' would allow Maker
to achieve a better score than the one stated in Theorem 2. This suggests
that Breaker should claim an element next to the element claimed by
Maker, and hence the initial board splits into two disjoint boards.
Hence it is natural to consider games $G$ that are disjoint union
of $F\left(l_{1}\right),\dots,F\left(l_{r}\right)$, $G\left(m_{1}\right),\dots,G\left(m_{s}\right)$
and $H\left(n_{1}\right),\dots,H\left(n_{t}\right)$. 

The plan of the paper is as follows. In Section 2 we prove a lower
bound for $\alpha\left(n,k\right)$ and deduce a lower bound for $\gamma_{b}\left(n\right)$.
In Section 3 we prove an upper bound for the score of games $G$ that
are disjoint union of games of the form $F\left(l\right)$, $G\left(m\right)$
and $H\left(n\right)$, and conclude the Theorems 1 and 2 from these
upper and lower bounds. 

\section{The lower bound}

Recall that $F\left(n,k\right)$ is defined to be the game played
on $\left\{ 1,\dots,n\right\} $, where at the start of the game Breaker
is allowed to claim $k$ elements, and then the players claim elements
in alternating order, and $\alpha\left(n,k\right)$ is the score attained
when both players play optimally. We start by proving the following
lower bound on $\alpha\left(n,k\right)$ which is later used to deduce
a lower bound on $\gamma_{b}\left(n\right)$. \\

\textbf{Lemma 3. $\alpha\left(n,k\right)\geq\left\lfloor \frac{n-3k+2}{5}\right\rfloor $}.
\\

\textbf{Proof. }Suppose that Breaker claims the elements $s_{1},\dots,s_{k}$
on her first move. These elements splits the path into $k+1$ (possibly
empty) intervals of lengths $l_{0},\dots,l_{k}$, with $l_{i}=s_{i+1}-s_{i}-1$
(with the convention $s_{0}=0$ and $s_{k+1}=n+1$). By symmetry we
may assume that $l_{0}$ is the longest interval. 

If $l_{0}\leq2$, then $n\leq k+2\cdot\left(k+1\right)=3k+2$, and
hence $\left\lfloor \frac{n-3k+2}{5}\right\rfloor =0$. Thus the claim
follows immediately in this case, and hence we may assume that $l_{0}\geq3$.
We treat the cases $l_{0}\geq4$ and $l_{0}=3$ individually. In both
cases the proof follows the same idea, however the choice of the initial
move is slightly different for $l_{0}=3$ since an interval with only
$3$ elements is 'too short' for the general argument. \\

\textbf{Case 1: $l_{0}\geq4$}. \\

The aim for Maker is to build a long block of consecutive elements
inside the interval. Initially, she claims the element $3$. Assuming
she has already claimed exactly the elements $\left\{ t,\dots,t+r\right\} $,
she claims one of $t+r+1$ or $t-1$, if possible. If not, she stops. 

Consider the point when this process terminates, and suppose that
at the point of termination she has claimed the set of elements $\left\{ t,\dots,t+r\right\} $.
Since one of these is the element $3$, we must have $t+r\geq3$ and
$t\in\left\{ 1,2,3\right\} $. Also note that the element $t+r+1$
must be claimed by Breaker, and also either $t=1$ or the element
$t-1$ is claimed by Breaker. Since $l_{0}\geq4$, it follows that
the elements $2$ and $4$ are not claimed after Maker's first move.
Since Breaker cannot claim both of these on her first move, it follows
that Maker can always guarantee that $r\geq1$. 

Let $T_{1}=\left\{ t+r+2,\dots n\right\} $ and let $b$ be the number
of elements claimed by Breaker in $T_{1}$. Note that Breaker has
claimed $k+r+1$ elements in total, and one of these must be $t+r+1$.
Furthermore, if $t>1$ then one of them must be $t-1$ as well. Hence
$b\leq k+r$, and if $t\geq2$ we also have $b\leq k+r-1$. Also note
that Maker has not claimed any elements in $T_{1}$. 

Note that claiming the elements $\left\{ t,\dots,t+r\right\} $ increases
the score by exactly $r$, and this is the only contribution for the
score coming outside $T_{1}$. Thus the total score that Maker can
attain is at least $r+\alpha\left(n-t-r-1,b\right)$. By induction,
it follows that the score is at least 
\begin{equation}
r+\left\lfloor \frac{n-t-r-1-3b+2}{5}\right\rfloor .\label{eq:0.5}
\end{equation}

If $t=1$, it follows that $b\leq k+r$. Also the condition $t+r\geq3$
implies that and $r\geq2$. Hence (\ref{eq:0.5}) implies that Maker
can guarantee that the score is at least 
\[
\left\lfloor \frac{n-3k+r+1-t}{5}\right\rfloor \geq\left\lfloor \frac{n-3k+2}{5}\right\rfloor 
\]
as required. 

If $t\geq2$, it follows that $b\leq k+r-1$. Recall that we always
have $t\le3$ and $r\geq1$. Hence (\ref{eq:0.5}) implies that Maker
can guarantee that the score is at least

\[
\left\lfloor \frac{n-3k-t+r+4}{5}\right\rfloor \geq\left\lfloor \frac{n-3k-3+1+4}{5}\right\rfloor =\left\lfloor \frac{n-3k+2}{5}\right\rfloor .
\]
Hence we have $\alpha\left(n,k\right)\geq\left\lfloor \frac{n-3k+2}{5}\right\rfloor $,
as required. $\square$\\

\textbf{Case 2: }$l_{0}=3$. \\

Again, Maker is aiming to claim as long block of consecutive elements
in $\left\{ 1,2,3\right\} $ as possible. Initially she claims the
element $2$. Since Breaker cannot pick both $1$ and $3$ on her
first move, Maker can always guarantee that the length of this block
is at least 2. If possible, she picks the last element on her third
move. 

Thus at the end of this process, exactly one of the following is true:
\begin{enumerate}
\item Maker has claimed all three elements in $\left\{ 1,2,3\right\} $.
\item Maker has claimed two consecutive elements in $\left\{ 1,2,3\right\} $
and Breaker has claimed the third element in $\left\{ 1,2,3\right\} $. 
\end{enumerate}
In both cases, consider the game played on $T_{1}=\left\{ 5,\dots,n\right\} $.
Let $a$ be the number of elements Maker claims in $\left\{ 1,2,3\right\} $.
Note that in both cases Breaker claims all the other elements in $\left\{ 1,2,3,4\right\} $
not claimed by Maker, and thus Breaker claims $4-a$ elements in $\left\{ 1,2,3,4\right\} $.
Since Breaker claims in total $a+k$ elements, it follows that she
claims $k+2a-4$ elements on $T_{1}$. Since Maker has not yet claimed
any elements in $T_{1}$, it follows that on $T_{1}$ Maker can increase
the score by $\alpha\left(n-4,k+2a-4\right)$. Since she has achieved
a score of $a-1$ outside $T_{1}$ with her block of $a$ consecutive
elements, it follows that the total score achieved is $a-1+\alpha\left(n-4,k+2a-4\right)$. 

By induction, it follows that the score achieved is at least

\[
a-1+\left\lfloor \frac{n-4-3\left(k+2a-4\right)+2}{5}\right\rfloor =\left\lfloor \frac{n-3k-a+5}{5}\right\rfloor .
\]
Since $a\in\left\{ 2,3\right\} $, it follows that 
\[
\alpha\left(n,k\right)\geq\left\lfloor \frac{n-3k+2}{5}\right\rfloor 
\]
as required. $\square$

Thus Lemma 3 holds by induction. $\Square$\\

\textbf{Lemma 4. $\gamma_{b}\left(n\right)\geq\left\lfloor \frac{n+4}{5}\right\rfloor $
}for $n\geq2$ and $\gamma_{b}\left(1\right)=0$. \\

\textbf{Proof. }When $n=1$, the claim is trivial as the only move
is given to Breaker. Now we consider the case $n\geq2$. 

At the start of the game, Maker is aiming to claim as long blocks
of consecutive elements as possible near the endpoints. Once this
is no longer possible, she starts using the same strategy as in Lemma
3. We start by describing this initial process formally. 

Suppose that after Maker's $k^{th}$ move the set of elements claimed
by Maker is of the form $\left\{ 1,\dots t\right\} \cup\left\{ n-k+t+1,\dots,n\right\} $
for some $t\in\left\{ 0,\dots,k\right\} $, with the convention that
$\left\{ 1,\dots,t\right\} =\emptyset$ when $t=0$ and $\left\{ n-k+t+1,\dots,n\right\} =\emptyset$
when $t=k$. Note that this certainly holds when $k=0$, as Maker
has not claimed any elements before her first move. If at least one
of the elements $t+1$ or $n-k+t$ is not yet claimed before Maker's
$k+1^{th}$ move, then Maker claims one of these elements which is
still available, and thus the set of vertices claimed by Maker is
of this form also after $k+1$ moves. If both $t+1$ and $n-k+t$
are claimed by Breaker, then the process stops.

This process terminates trivially, as Breaker must claim an element
during the game. Suppose that when the process terminates, the set
of vertices claimed by Maker is of the form $\left\{ 1,\dots t\right\} \cup\left\{ n-k+t+1,\dots,n\right\} $
for some $k$ and $t\in\left\{ 0,\dots,k\right\} $. Note that we
must have $k\geq1$, as Breaker cannot claim both elements $1$ and
$n$ on her first move. 

Let $T=\left\{ t+2,\dots,n-k+t-1\right\} $, and note that by the
choice of $k$ and $t$ it follows that Maker has not claimed any
elements in $T$. Since the process has terminated at this stage,
it follows that Breaker must have claimed the elements $t+1$ and
$n-k+t$. Since Breaker started the game, she has claimed $k+1$ elements
in total, and thus $k-1$ of these elements must be in $T$. 

Note that any increment of the score arising outside $T$ occurs from
the sets $\left\{ 1,\dots,t\right\} $ and $\left\{ n-k+t+1,\dots,n\right\} $.
On the other hand, since Maker has not claimed any elements in $T$
and Breaker has claimed $k-1$ elements in $T$, the rest of the game
on $T$ corresponds to the game $F\left(n-k-2,\,k-1\right)$. Hence
Maker can increase the score by at least $\alpha\left(n-k-2,k-1\right)$
in $T$.

It is easy to check that the contribution on the score arising from
the intervals $\left\{ 1,\dots,t\right\} $ and $\left\{ n-k+t+1,\dots,n\right\} $
is exactly $t+\left(k-t\right)=k$. Hence by Lemma 3, it follows that
Maker can guarantee that the score is at least
\[
k+\alpha\left(n-k-2,k-1\right)\geq k+\left\lfloor \frac{n-k-2-3\left(k-1\right)+2}{5}\right\rfloor =\left\lfloor \frac{n+k+3}{5}\right\rfloor .
\]
Since $k\geq1$, it follows that Maker can always guarantee that the
score is at least $\left\lfloor \frac{n+4}{5}\right\rfloor $, which
completes the proof. $\square$

\section{The upper bound }

In this section, all congruences are considered modulo $5$ unless
otherwise stated, and in such cases we omit $(\text{mod 5})$ from
the notation. Furthermore, we write $n\equiv0\text{ or }1$ instead
of $n\equiv0$ or $n\equiv1$, and $n\not\equiv0\text{ and }1$ instead
of $n\not\equiv0$ and $n\not\equiv1$.\\

\textbf{Lemma 5. }Suppose $T$ is a disjoint union of games $F\left(l_{1}\right),\dots,F\left(l_{r}\right)$,
$G\left(m_{1}\right),\dots,G\left(m_{s}\right)$ and $H\left(n_{1}\right),\dots,H\left(n_{t}\right)$,
with Maker having the first move. Let $f\left(\underline{l};\,\underline{m};\,\underline{n}\right)$
be the score of this game when both players play optimally. Let $N_{1}=\left|\left\{ i\,:\,l_{i}\equiv3\text{ or }4\right\} \right|$,
$N_{2}=\left|\left\{ i\,:\,m_{i}\equiv0\text{ or }1\right\} \right|$,
$N_{3}=\left|\left\{ i\,:n_{i}\neq2\text{ and }\,n_{i}\equiv2\text{ or }3\right\} \right|$,
$N_{4}=\left|\left\{ i\,:\,n_{i}=2\right\} \right|$ and $N_{5}=\left|\left\{ i\,:\,n_{i}=1\right\} \right|$.
Let $\epsilon\in\left\{ 0,1\right\} $ be chosen such that $N_{5}\equiv\epsilon\text{\text{ (\text{mod 2})}}$.
Then we have 
\begin{equation}
f\left(\underline{l};\,\underline{m};\,\underline{n}\right)\le\sum_{i=1}^{r}\left\lfloor \frac{l_{i}+2}{5}\right\rfloor +\sum_{i=1}^{s}\left\lfloor \frac{m_{i}+5}{5}\right\rfloor +\sum_{i=1}^{t}\left\lfloor \frac{n_{i}+8}{5}\right\rfloor -N_{4}+\epsilon-\left\lfloor \frac{N_{1}+N_{2}+N_{3}+\epsilon}{2}\right\rfloor .\label{eq:1}
\end{equation}
\\

By looking at the proof of Lemma 3, it is reasonable for Breaker to
claim one of the points next to the point Maker claimed on her first
move, as in this case Breaker can restrict the length of intervals
created by Maker. Such a first pair of moves splits the original board
into two new boards, which motivates the idea of considering unions
of disjoint boards. It might be tempting to say, that Breaker can
always follow Maker into the board where she plays her next move,
and hence proceed by using an inductive proof. However, sometimes
Breaker may gain an 'extra move' if one of these boards has no sensible
moves left (i.e.\ the component is $F\left(1\right)$ or $F\left(2\right)$). 

Ignoring these extra moves completely would make the proof much shorter,
but the bound obtained that way would not even be good enough asymptotically.
Since Maker is free to alternate between these two boards, she has
some control on the time of the game when Breaker is given this extra
move. In particular, in this case we cannot assume that these extra
moves are given at the start of the game, which was the case in Section
2. In order to keep track of these extra moves, we need to consider
arbitrary disjoint unions of boards. 

We start by briefly outlining the structure of the proof and explaining
where the upper bound in (\ref{eq:1}) comes from. The proof is by
induction on the sum of the lengths of the paths. The aim is to prove
that for any possible Maker's initial move, there is a move for Breaker
that can be used to show that (\ref{eq:1}) holds by induction. This
move will in general depend on the position of the initial move modulo
$5$, however we have to be slightly more careful if the initial move
is close to the endpoints of a board. For the same reason, one has
to be careful with small components of the board as well. 

Since there are $3$ possible board types, $5$ possible locations
for the initial move (mod 5), and two possible cases for the size
of the initial length of the component (depending on whether the initial
length is involved in one of $N_{1}$, $N_{2}$ or $N_{3}$, or not),
it follows that there are in some sense $30$ cases to be considered.
In addition, we have to cover small cases as well. Fortunately, some
of these cases can be treated simultaneously, and in general the techniques
used to prove various cases are identical or use very similar techniques.

In a sense the hardest part is rather to come up with a suitable upper
bound in (\ref{eq:1}) that is strong enough for an inductive argument
to work than the proof itself. Once a suitable upper bound is chosen,
identifying possible 'response moves' for Breaker is reasonably easy.
Finally, the proof itself is mathematically not challenging, but it
is reasonably tedious. 

Why should we choose this particular upper bound in (\ref{eq:1})?
For $B=F\left(l\right)$, $G\left(m\right)$ or $H\left(n\right)$
(with $n\geq3$) it turns out that Breaker can always guarantee that
the score is at most $\left\lfloor \frac{l+2}{5}\right\rfloor $,
$\left\lfloor \frac{m+5}{5}\right\rfloor $ or $\left\lfloor \frac{n+8}{5}\right\rfloor $
respectively. This explains the first three sums in the upper bound.
Moreover, if $l\equiv3\text{ or }4$, $m\equiv0\text{ or }1$ or $n\equiv2\text{ or }3$
it turns out that Breaker has a strategy which allows her to force
Maker to either play the last non-trivial move (i.e.\ after which all
components are either empty, $F\left(1\right)$ or $F\left(2\right)$),
or Maker can only attain a score which is strictly less than this
bound. Hence the quantity $N_{1}+N_{2}+N_{3}$ is measuring the number
of these 'additional moves'. Given such an additional move Breaker
can make another component of the board slightly shorter, which either
reduces the score by one or guarantees that she will gain an extra
move from that board as well. 

However, one has to be careful with small values of $n$. Indeed,
it turns out that on $H\left(2\right)$ Maker can only increase the
score by $1$ (instead of $2$), and Breaker cannot gain an extra
turn. This is the reason behind the $-N_{4}$-term. Also on $H\left(1\right)$
Maker can score $2$ points (instead of $1$), and Breaker gains an
extra turn. Note that if the number of components of the form $H\left(1\right)$
is even, then Breaker can always claim a point on another component
that is $H\left(1\right)$. If the number is odd, she can follow this
pairing strategy until the number of such boards decreases to $1$,
in which case she has to use the extra move elsewhere. This is the
reason behind the fact that only the parity of $N_{5}$ matters. 

In a sense, dealing with boards of the form $H\left(n\right)$ is
the hardest task due to irregular behaviour of these boards when $n$
is small. Hence we start the proof by considering these type of boards,
and during the proof we also introduce some standard arguments that
can be easily used when dealing with boards of the form $F\left(l\right)$
or $G\left(m\right)$. In those cases, we do not always give full
justification.

Note that the bound (\ref{eq:1}) may not always be tight, but by
a similar argument as presented in Section 2 one could verify that
it is tight when applied to $F\left(l\right)$, $G\left(m\right)$
or $F\left(n\right)$, which is good enough for our purposes. The
reason why the bound is not necessarily tight is the fact that sometimes
Breaker could have a better place to play her extra move, rather than
the 'worst case scenario' that is considered in the proof.

For convenience define 
\[
g\left(\underline{l};\,\underline{m};\,\underline{n}\right)=\sum_{i=1}^{r}\left\lfloor \frac{l_{i}+2}{5}\right\rfloor +\sum_{i=1}^{s}\left\lfloor \frac{m_{i}+5}{5}\right\rfloor +\sum_{i=1}^{t}\left\lfloor \frac{n_{i}+8}{5}\right\rfloor -N_{4}+\epsilon-\left\lfloor \frac{N_{1}+N_{2}+N_{3}+\epsilon}{2}\right\rfloor ,
\]
 
\[
y\left(\underline{l};\,\underline{m};\,\underline{n}\right)=\sum_{i=1}^{r}\left\lfloor \frac{l_{i}+2}{5}\right\rfloor +\sum_{i=1}^{s}\left\lfloor \frac{m_{i}+5}{5}\right\rfloor +\sum_{i=1}^{t}\left\lfloor \frac{n_{i}+8}{5}\right\rfloor 
\]
 and 
\[
z\left(\underline{l};\,\underline{m};\,\underline{n}\right)=-N_{4}+\epsilon-\left\lfloor \frac{N_{1}+N_{2}+N_{3}+\epsilon}{2}\right\rfloor .
\]
For later purposes, it is convenient to observe that we may rewrite
$z$ as 
\begin{equation}
z\left(\underline{l};\,\underline{m};\,\underline{n}\right)=-N_{4}-\left\lfloor \frac{N_{1}+N_{2}+N_{3}-\epsilon}{2}\right\rfloor .\label{eq:3}
\end{equation}
\\

\textbf{Proof. }Define $N=\sum_{i=1}^{r}l_{i}+\sum_{i=1}^{s}m_{i}+\sum_{i=1}^{t}n_{i}$.
The proof is by induction on $N$, and it is easy to check that the
claim holds for all possible configurations when $N=1$ or $N=2$.
Suppose the claim holds whenever $N\leq M-1$ for some $M\geq3$,
and suppose that $\underline{l},\,\underline{m},\,\underline{n}$
are chosen such that $\sum_{i=1}^{r}l_{i}+\sum_{i=1}^{s}m_{i}+\sum_{i=1}^{t}n_{i}=M$. 

We now split the proof into several cases depending on Maker's first
move. In each case, let $S\left(T\right)$ be the maximum score that
Maker can attain given this first move and given that Breaker plays
optimally. \\

\textbf{Case 1: }Maker plays on on $H\left(n_{t}\right)$.\\

For convenience set $n=n_{t}$. The game $H\left(n\right)$ is played
on $\left\{ 1,\dots,n\right\} $, and since both endpoints of the
board are symmetric we may assume that Maker claims first an element
$j$ satisfying $j\leq\left\lceil \frac{n}{2}\right\rceil $. We prove
that apart from small values of $n$, claiming one of $j-1$ or $j+1$
is a suitable choice for Breaker, where the choice is made depending
on $j\text{ (mod 5)}$, as indicated in Table 1. If $j\geq3$, after
such first pair of moves it is easy to see that the $H\left(n\right)$-component
of the board splits into disjoint union of $H\left(a\right)$ and
$G\left(b\right)$ for some $a,\,b$ with $n=a+b+2$. However, since
the boards $H\left(1\right)$ and $H\left(2\right)$ behave in a different
way compared to other boards of the form $H\left(n\right)$, it turns
out to be convenient to consider the cases $j=1$, $j=2$ and $\left(j,n\right)=\left(3,5\right)$
individually. 

Indeed, if $4\leq j\leq\left\lceil \frac{n}{2}\right\rceil $, then
the board splits into $H\left(a\right)$ and $G\left(b\right)$ with
$a\geq3$. If $j=3$, then as indicated in Table 1 Breaker claims
the element $2$. Hence the boards splits into $G\left(1\right)$
and $H\left(n-3\right)$, which is one of $H\left(1\right)$ or $H\left(2\right)$
only if $n=5$, as $j\leq\left\lceil \frac{n}{2}\right\rceil $. Hence
$j=1$, $j=2$ and $\left(j,n\right)=\left(3,5\right)$ are the only
special cases which could change the number of boards of the form
$H\left(1\right)$ or $H\left(2\right)$. 

Denote the new set of parameters obtained after the first pair of
moves as $\underline{l}'$, $\underline{m}'$ and $\underline{n}'$,
and let $s_{i}$ denote the increment of the score caused by Maker's
first move. Throughout the proof it is convenient to define the quantities
$d_{1}=z\left(\underline{l};\,\underline{m};\,\underline{n}\right)-z\left(\underline{l}';\,\underline{m}';\,\underline{n}'\right)$
and $d_{2}=y\left(\underline{l};\,\underline{m};\,\underline{n}\right)-y\left(\underline{l}';\,\underline{m}';\,\underline{n}'\right)$.
Note that $g\left(\underline{l};\,\underline{m};\,\underline{n}\right)=d_{1}+d_{2}+g\left(\underline{l}';\,\underline{m}';\,\underline{n}'\right)$. 

By induction we know that $S\left(T\right)\leq g\left(\underline{l}';\,\underline{m}';\,\underline{n}'\right)+s_{i}$.
Since our aim is to prove that $S\left(T\right)\leq g\left(\underline{l};\,\underline{m};\,\underline{n}\right)$,
it suffices to prove that we always have $d_{1}+d_{2}\geq s_{i}$.
In fact, we will prove that for all possible Maker's initial moves
there exists a move for Breaker that satisfies $d_{1}+d_{2}\geq s_{i}$. 

We start with the general case $j\geq3$ and $n\geq6$, and we deal
with the special cases later. 

\begin{table}[h]
\caption{Choices for Breaker's first move depending on $j$}

~\\

\begin{tabular}{|c|c|c|c|c|c|c|}
\hline 
 & $F\left(n\right)$ & Condition on $a$ or $b$ & $G\left(n\right)$ & Condition on $a$ or $c$ & $H\left(n\right)$ & Condition on $a$ or $b$\tabularnewline
\hline 
$j\equiv0$ & $j+1$ & $b\equiv4$ & $j-1$ & $a\equiv3$ & $j-1$ & $b\equiv3$\tabularnewline
\hline 
$j\equiv1$ & $j-1$ & $a\equiv4$ & $j-1$ & $a\equiv4$ & $j-1$ & $b\equiv4$\tabularnewline
\hline 
$j\equiv2$ & $j+1$ & $b\equiv1$ & $j+1$ & $c\equiv1$ & $j+1$ & $a\equiv1$\tabularnewline
\hline 
$j\equiv3$ & $j-1$ & $a\equiv1$ & $j-1$ & $a\equiv1$ & $j-1$ & $b\equiv1$\tabularnewline
\hline 
$j\equiv4$ & $j-1$ & $a\equiv2$ & $j+1$ & $c\equiv3$ & $j+1$ & $a\equiv3$\tabularnewline
\hline 
\end{tabular}
\end{table}

\textbf{Case 1.1: }$n\geq6$, $j\geq3$.\\

In this case we have $s_{i}=0$, so it suffices to prove that $d_{1}+d_{2}\geq0$.
It is easy to see that $N_{1}$, $N_{4}$ and $\epsilon$ are unaffected
in this case. Since $N_{2}$ certainly cannot increase and $N_{3}$
can decrease by at most $1$, it follows that $d_{1}\geq-1$. 

Note that we have $d_{2}=\left\lfloor \frac{n+8}{5}\right\rfloor -\left\lfloor \frac{a+8}{5}\right\rfloor -\left\lfloor \frac{b+5}{5}\right\rfloor $.
By using the trivial upper and lower bounds $x-1\leq\left\lfloor x\right\rfloor \leq x$
and the fact that $n=a+b+2$, it follows that $d_{2}\geq\frac{n+3}{5}-\frac{a+b+13}{5}=\frac{-8}{5}$.
Since $d_{2}$ is an integer, it follows that $d_{2}\geq-1$. We now
split to several two cases based on the value of $n\text{ (mod 5)}$
in order to improve our bounds on $d_{1}$ and $d_{2}$ to attain
$d_{1}+d_{2}\geq0$. \\

\textbf{Case 1.1.1: }$n\equiv2\text{ or }3$.\\

We start by improving the bound on $d_{2}$. Since $n\equiv2\text{ or }3$
it follows that $\left\lfloor \frac{n+8}{5}\right\rfloor \geq\frac{n+7}{5}$.
Hence by using the trivial bounds for the other terms, we obtain that
$d_{2}\geq\frac{n+7}{5}-\frac{a+b+13}{5}=\frac{-4}{5}$. Since $d_{2}$
is an integer, it follows that $d_{2}\geq0$. 

First suppose that $a\equiv2\text{ or }3$. Then $N_{3}$ cannot decrease,
so in fact we have $d_{1}\geq0$. Hence we have $d_{1}+d_{2}\geq0$,
as required. 

Now suppose that $b\equiv0\text{ or }1$. Then $N_{3}$ decreases
by at most $1$ and $N_{2}$ increases by $1$. Hence the sum $N_{2}+N_{3}$
certainly cannot decrease. Thus we also have $d_{1}\geq0$, and thus
it follows that $d_{1}+d_{2}\geq0$, as required. 

Finally suppose that $a\not\equiv2\text{ and }3$ and $b\not\equiv0\text{ and }1$.
Then we have $\left\lfloor \frac{a+8}{5}\right\rfloor +\left\lfloor \frac{b+5}{5}\right\rfloor \leq\frac{a+6}{5}+\frac{b+3}{5}=\frac{a+b+9}{5}$.
Note that the equality holds if and only if $a\equiv4$ and $b\equiv2$,
but by Table 1 it follows that this can never happen. Hence this inequality
must be strict, and hence it follows that $d_{2}>\frac{n+7}{5}-\frac{a+b+9}{5}=0$.
Hence we must have $d_{2}\geq1$, and combining this with the trivial
bound $d_{1}\geq-1$ it follows that $d_{1}+d_{2}\geq0$, as required.
This completes the proof of Case 1.1.1. \\

\textbf{Case 1.1.2: }$n\not\equiv2\text{ and }3$.\\

Since $n\not\equiv2\text{ and }3$, it follows that $N_{3}$ cannot
decrease. Hence we must have $d_{1}\geq0$. 

First suppose that $a\equiv2\text{ or }3$ and $b\equiv0\text{ or }1$.
Then both $N_{2}$ and $N_{3}$ increase by $1$, and hence it follows
that $d_{1}\geq1$. Combining this with the trivial bound $d_{2}\geq-1$
implies that $d_{1}+d_{2}\geq0$, as required.

Now suppose that $a\not\equiv2\text{ and }3$ or $b\not\equiv0\text{ and }1$.
As in Case 1.1.1, in both cases we can improve the upper bound on
$\left\lfloor \frac{a+8}{5}\right\rfloor +\left\lfloor \frac{b+5}{5}\right\rfloor $
to $\left\lfloor \frac{a+8}{5}\right\rfloor +\left\lfloor \frac{b+5}{5}\right\rfloor \leq\frac{a+b+11}{5}$,
and note that the equality holds if and only if ($a\equiv4$ and $b\equiv0$)
or ($a\equiv2$ and $b\equiv2$). However, note that by Table 1 both
of these cases are impossible. Hence the inequality must be strict,
and thus we have $d_{2}>\frac{n+4}{5}-\frac{a+b+11}{5}=-1$. Hence
it follows that $d_{2}\geq0$, and thus we have $d_{1}+d_{2}\geq0$,
which completes the proof of Case 1.1.2. \\

\textbf{Case 1.2: }$j=1$.\\

Here we split into three cases based on the size of $n$. First, we
consider the case $n\geq3$ which should be viewed as the main part
of the argument. Then we consider the cases $n=2$ and $n=1$ individually,
as these behave in a slightly different way as the boards are small.
The case $n=1$ turns out to be very tedious and lengthy, and it does
not really contain any interesting ideas either. In some sense, the
only task in this case is to find out a good enough way for Breaker
to use her additional move.\\

\textbf{Case 1.2.1: }$n\geq3$. \\

Suppose Breaker claims the element $2$. Since Maker claimed $1$
on board $H\left(n\right)$ with $n\geq3$, it follows that $s_{i}=1$.
Hence it suffices to prove that with this move Breaker can achieve
$d_{1}+d_{2}\geq1$. First of all, note that the board $H\left(n\right)$
was replaced by $G\left(n-2\right)$, which is non-empty as $n\geq3$.
Since $n\geq3$, and $n\equiv2\text{ or }3$ if and only if $n-2\equiv0\text{ or }1$,
it follows that $N_{3}$ decreases by $1$ if and only if $N_{2}$
increases by $1$. In particular, it follows that $d_{1}=0$. On the
other hand, it is easy to check that $d_{2}=\left\lfloor \frac{n+8}{5}\right\rfloor -\left\lfloor \frac{\left(n-2\right)+5}{5}\right\rfloor =1$.
Hence we always have $d_{1}+d_{2}=1$, which completes the proof of
Case 1.2.1.\\

\textbf{Case 1.2.2: }$n=2$.\\

Suppose Breaker claims the element $2$. Since the board $H\left(2\right)$
has only two elements, it follows that all elements of the board are
occupied after this pair of moves. Note that we certainly have $s_{i}=1$,
and $d_{2}=\left\lfloor \frac{2+8}{5}\right\rfloor =2$. On the other
hand, it is clear that $N_{1}$, $N_{2}$, $N_{3}$ and $\epsilon$
remain unaffected while $N_{4}$ decreases by $1$. Hence we have
$d_{1}=-1$, and thus $d_{1}+d_{2}=1$ as required. This completes
the proof of Case 1.2.2. \\

\textbf{Case 1.2.3: }$n=1$. \\

Since $n=1$, it follows that $s_{i}=2$. First suppose $N_{5}>1$,
and that Breaker chooses another board of the form $H\left(1\right)$
and claims the only element on that board. Hence $N_{5}$ decreases
by $2$, so $\epsilon$ remains unaffected and thus we have $d_{1}=0$.
On the other hand, we have $d_{2}=2\left\lfloor \frac{1+8}{5}\right\rfloor =2$,
and hence it follows that $d_{1}+d_{2}=s_{i}$, as required. 

Otherwise we must have $N_{5}=1$, and hence we certainly have $\epsilon=1$.
Since the total number of points on $T$ is strictly more than $1$,
it follows that there exists another component $B$ of $T$. 

First suppose that $B=H\left(2\right)$ and that Breaker claims the
element $1$. Then $N_{2}$ increases by $1$, $N_{4}$ decreases
by $1$ and $\epsilon$ is replaced by $1-\epsilon$. Hence $N_{1}+N_{2}+N_{3}-\epsilon$
increases by $2$ and $-N_{4}$ increases by $1$, so we have $d_{1}=0$.
Note that $d_{2}=\left\lfloor \frac{1+8}{5}\right\rfloor +\left\lfloor \frac{2+8}{5}\right\rfloor -\left\lfloor \frac{1+5}{5}\right\rfloor =2$,
and thus it follows that $d_{1}+d_{2}=s_{i}$, as required. 

Now suppose that $B=H\left(m\right)$ with $m\geq3$. Suppose that
Breaker claims the element $1$, and hence $B$ is replaced by $G\left(m-1\right)$.
Then $N_{4}$ remains unaffected, $N_{3}$ decreases by at most $1$
and $N_{2}$ increases by at most one. Since $\epsilon$ changes from
$1$ to $0$, it follows that $N_{1}+N_{2}+N_{3}-\epsilon$ cannot
decrease, and hence we have $d_{1}\geq0$. Note that we have $d_{2}=\left\lfloor \frac{1+8}{5}\right\rfloor +\left\lfloor \frac{m+8}{5}\right\rfloor -\left\lfloor \frac{m+4}{5}\right\rfloor $,
and thus we trivially have $d_{2}\geq1$. 

If $m\equiv1$, then $m-1\equiv0$ and thus $N_{2}$ increases by
$1$ but $N_{3}$ does not decrease. Hence $N_{1}+N_{2}+N_{3}-\epsilon$
increases by $2$, and thus $d_{1}\geq1$. If $m\not\equiv1$, then
we certainly have $d_{2}\geq2$. Hence in either case we have $d_{1}+d_{2}\geq2$,
as required. 

Next suppose that $B=G\left(m\right)$, and suppose that Breaker claims
the element $1$. Hence $B$ is replaced by $F\left(m-1\right)$.
As above, it is easy to deduce that $N_{4}$ remains unaffected and
$N_{1}+N_{2}+N_{3}-\epsilon$ cannot decrease, and hence we have $d_{1}\geq0$.
We also have $d_{2}=\left\lfloor \frac{1+8}{5}\right\rfloor +\left\lfloor \frac{m+5}{5}\right\rfloor -\left\lfloor \frac{m+1}{5}\right\rfloor $,
and thus $d_{2}\geq1$. 

If $m\equiv4$, then $m-1\equiv3$ and thus $N_{1}$ increases by
$1$ but $N_{2}$ does not decrease. Hence we can similarly deduce
that $d_{1}\geq1$. Otherwise, it is easy to see that $d_{2}\geq2$.
Hence in either case we have $d_{1}+d_{2}\geq2$, as required. Note
that the same argument also applies even when $m=1$ (with the convention
that $F\left(0\right)$ is an empty board). 

Finally suppose that $B=F\left(m\right)$, and suppose that Breaker
claims the element $m-2$. Hence $B$ is replaced by disjoint union
of $F\left(m-3\right)$ and $F\left(2\right)$, but the component
of the form $F\left(2\right)$ can be omitted as on this board Breaker
can follow pairing strategy to avoid any increment in the score. Again,
we know that $N_{1}$ cannot decrease by more than $1$, and hence
$d_{1}\geq0$. We also have $d_{2}=\left\lfloor \frac{1+8}{5}\right\rfloor +\left\lfloor \frac{m+2}{5}\right\rfloor -\left\lfloor \frac{m-1}{5}\right\rfloor $,
and hence $d_{2}\geq1$. 

If $m\equiv3,\,4\text{ or }5$, then we certainly have $d_{2}\geq2$.
If $m\equiv1\text{ or }2$, then $m-3\equiv3\text{ or }4$, and hence
$N_{1}$ increases by $1$. Hence $N_{1}+N_{2}+N_{3}-\epsilon$ increases
by 2, and thus we must have $d_{1}\geq1$. Hence in either case we
have $d_{1}+d_{2}\geq2$. This completes the proof of Case 1.2.3.
\\

\textbf{Case 1.3: $j=2$}.\\

Since $j\le\left\lceil \frac{n}{2}\right\rceil $, it follows that
we must have $n\geq3$. Hence we split into cases based on whether
$n\geq5$, $n=4$ or $n=3$. \\

\textbf{Case 1.3.1: $n\geq5$}. \\

Suppose Breaker claims the element $1$. Hence $s_{i}=0$, and the
board becomes a copy of $H\left(n-2\right)$. Since $n-2\geq3$, it
follows that $N_{4}$ and $\epsilon$ remain unchanged. 

If $n\not\equiv2\text{ and }3$, then $N_{3}$ cannot decrease. Hence
it follows that $d_{1}\geq0$. We also have $d_{2}=\left\lfloor \frac{n+8}{5}\right\rfloor -\left\lfloor \frac{n+6}{5}\right\rfloor \geq0$,
and thus $d_{1}+d_{2}\geq0$ as required. 

If $n\equiv2\text{ or }3$, then $N_{3}$ decreases by at most $1$
and hence we have $d_{1}\geq-1$. We again have $d_{2}=\left\lfloor \frac{n+8}{5}\right\rfloor -\left\lfloor \frac{n+6}{5}\right\rfloor $,
and since $n\equiv2\text{ or }3$ it follows that $d_{2}\geq1$. Thus
$d_{1}+d_{2}\ge0$, which completes the proof of Case 1.3.1. \\

\textbf{Case 1.3.2: }$n=4$. \\

Again suppose that Breaker claims the element $1$. Hence $s_{i}=0$,
and since $4\not\equiv2\text{ and }3$ it follows that $N_{1}$, $N_{2}$,
$N_{3}$ and $\epsilon$ remain unaffected. On the other hand, by
definition we know that $N_{4}$ increases by $1$ as after this pair
of moves the board becomes $H\left(2\right)$. Hence we have $d_{1}=1$.
We also have $d_{2}=\left\lfloor \frac{4+8}{5}\right\rfloor -\left\lfloor \frac{2+8}{5}\right\rfloor =0$,
and thus it follows that $d_{1}+d_{2}=1>0$, which completes the proof
of Case 1.3.2. \\

\textbf{Case 1.3.3: $n=3$}. \\

Again suppose that Breaker claims the element $1$, and thus we have
$s_{i}=0$. Note that after this pair of moves we are left with $H\left(1\right)$,
and it is easy to verify that $d_{2}=\left\lfloor \frac{3+8}{5}\right\rfloor -\left\lfloor \frac{1+8}{5}\right\rfloor =1$. 

It is clear that $N_{1}$, $N_{2}$ and $N_{4}$ remain unchanged.
It is easy to observe that $N_{3}$ decreases by $1$, and $\epsilon$
is replaced by $1-\epsilon$. Hence in the worst case $N_{3}-\epsilon$
decreases by $2$, and thus by (\ref{eq:3}) it follows that $d_{1}\geq-1$,
and hence we have $d_{1}+d_{2}\geq0$. This completes the proof of
Case 1.3.3. \\

\textbf{Case 1.4: }$n=5$ and $j=3$. \\

Suppose that Breaker claims the element $2$. Hence the board $B\left(5\right)$
splits into $G\left(1\right)$ and $H\left(2\right)$, and we have
$s_{i}=0$. Hence $N_{4}$ increases by $1$, $N_{2}$ increases by
$1$ and $N_{1}$, $N_{3}$ and $\epsilon$ remain constant. Thus
we must have $d_{1}\geq1$. On the other hand, note that $d_{2}=\left\lfloor \frac{5+8}{5}\right\rfloor -\left\lfloor \frac{2+8}{5}\right\rfloor -\left\lfloor \frac{1+5}{5}\right\rfloor =-1$.
Hence it follows that $d_{1}+d_{2}\geq0$, which completes the proof
of Case 1.4. \\

This completes the proof of Case 1. $\square$\\

\textbf{Case 2:} Maker plays on on $G\left(m_{s}\right)$. \\

For convenience set $n=m_{s}$. The game $G\left(n\right)$ is played
on $\left\{ 1,\dots,n\right\} $, and note that in this case the board
is not symmetric. Hence we choose the labeling so that claiming the
element $1$ increases the score by $1$, but claiming the element
$n$ does not. 

Assume that Maker plays her first move in position $j$. As before
we prove that claiming $j-1$ or $j+1$ is a suitable choice for Breaker,
and this choice is again determined by $j$ $\text{(mod 5)}$. We
use the same notation as before, however in this cases there are two
options on how the board might split: the board either splits into
components of the form $G\left(a\right)$ and $G\left(b\right)$ if
Breaker claims $j-1$, or into components of the form $H\left(c\right)$
and $F\left(d\right)$ if Breaker claims $j+1$. In this case we only
need to consider the cases $j=1$, $j=2$ and $j=n$ individually,
and note that hence we may assume that $n\geq4$. We start by checking
the special cases, and we skip some of the details when they are identical
to the arguments used in Case 1. \\

\textbf{Case 2.1: $j=1$}\\
. 

This is essentially identical to the proof of Case 1.2.1. Indeed,
suppose Breaker claims the element 2. Hence after the first pair of
moves the board becomes $F\left(n-2\right)$ and we have $s_{i}=1$.
As in the proof of Case 1.2.1, we have $s_{i}=1$ and $d_{2}=\left\lfloor \frac{n+5}{5}\right\rfloor -\left\lfloor \frac{n}{5}\right\rfloor =1$.
Also as in Case 1.2.1 it follows that $N_{2}$ decreases by $1$ if
and only if $N_{1}$ increases by $1$, and hence $d_{1}=0$. Thus
$d_{1}+d_{2}=1$, as required. \\

\textbf{Case 2.2: $j=2$}.\\

This is identical to the proof of Case 1.3.1. \\

\textbf{Case 2.3: }$j=n$. \\

Suppose that Breaker claims the element  $n-1$. After this pair of
moves the board becomes $G\left(n-2\right)$ and we have $s_{i}=0$.
Note that $N_{2}$ can decrease by at most $1$, and hence $d_{1}\geq-1$.
We also have $d_{2}=\left\lfloor \frac{n+5}{5}\right\rfloor -\left\lfloor \frac{n+3}{5}\right\rfloor $,
and thus we certainly have $d_{2}\geq0$. 

If $n\equiv0\text{ or }1$, it is easy to verify that we have $d_{2}=1$,
and hence $d_{1}+d_{2}\geq0$ as required. Otherwise it follows  that
$N_{2}$ cannot decrease, and hence we must have $d_{1}\geq0$. Thus
$d_{1}+d_{2}\geq0$ holds in this case as well, which completes the
proof of Case 2.3. \\

\textbf{Case 2.4: }$3\leq j\leq n-1$. \\

Suppose that Breaker chooses the appropriate move indicated in Table
1 depending on the value of $j\text{ (mod 5)}$. Note that depending
on the value of $j$, the board may split into components of the form
$G\left(a\right)$ and $G\left(b\right)$ or of the form $H\left(c\right)$
and $F\left(d\right)$. We now consider 4 cases, depending on the
value value of $n\text{ (mod 5)}$ and depending on how the board
splits into two components. As in Case 1.1, we have the trivial lower
bounds $d_{1}\geq-1$ and $d_{2}\geq-1$. \\

\textbf{Case 2.4.1:} $n\equiv0\text{ or }1$.\\

Note that regardless of how the board splits into two components,
we can deduce in either case by using the trivial upper bound $\left\lfloor x\right\rfloor \leq x$
that $d_{2}\geq\left\lfloor \frac{n+5}{5}\right\rfloor -\frac{n+8}{5}$.
Since $n\equiv0\text{ or }1$, it follows that $d_{2}\geq\frac{-4}{5}$
and thus $d_{2}\geq0$. \\

\textbf{Case 2.4.1.1:} $j\equiv0,\,1\text{ or }3$.\\

In this case Breaker claims the element $j-1$, and hence the board
splits into components of the form $G\left(a\right)$ and $G\left(b\right)$.
Note that from Table 1 we can conclude that $a\equiv1,\,3\text{ or }4$. 

First suppose that $a\equiv1$ or $b\equiv0\text{ or }1$. Then $N_{2}$
certainly does not decrease, so $d_{1}\geq0$. Hence $d_{1}+d_{2}\geq0$,
as required. 

Otherwise we must have $a\equiv3\text{ or }4$ and $b\equiv2,\,3\text{ or }4$.
Hence we must have $\left\lfloor \frac{a+5}{5}\right\rfloor +\left\lfloor \frac{b+5}{5}\right\rfloor \leq\frac{a+2}{5}+\frac{b+3}{5}=\frac{n+3}{5}$.
Since $n\equiv0\text{ or }1$, it follows that $d_{2}\geq\frac{n+4}{5}-\frac{n+3}{5}>0$,
and thus $d_{2}\geq1$. Hence $d_{1}+d_{2}\geq0$, which completes
the proof Case 2.4.1.1.\\

\textbf{Case 2.4.1.2:} $j\equiv2\text{ or }4$.\\

In this case Breaker claims the element $j+1$ and the board splits
into components of the form $H\left(c\right)$ and $F\left(d\right)$.
Since $j>2$, it follows that $j\geq4$ and thus $c\geq3$. Hence
$N_{4}$ and $\epsilon$ are unaffected by the first pair of moves.
Again, we will split into cases depending on whether one of $c\equiv2\text{ or }3$
or $d\equiv3\text{ or }4$ holds or not. The details follow exactly
as in Case 2.4.1.1, and hence we omit the proof. \\

\textbf{Case 2.4.2:} $n\not\equiv0\text{ and }1$. \\

Now regardless of how the board splits into two components we can
deduce that $d_{1}\geq0$, as none of the $N_{i}$'s can decrease.
Again, the rest of the proof is similar to the proof of Case 1.1.2
(with appropriate modifications similar to those done in Case 2.4.1).
Hence we skip the details. \\

This completes the proof of Case 2. $\square$\\

\textbf{Case 3: }Maker plays on on $F\left(l_{r}\right)$. \\

For convenience set $n=l_{r}$. The game $F\left(n\right)$ is played
on $\left\{ 1,\dots,n\right\} $, and this time the board is again
symmetric. Hence we may assume that Maker plays her first move $j$
in a position with $j\leq\left\lceil \frac{n}{2}\right\rceil $.This
time the only special case that needs to be considered is $j=1$,
and again we prove that for $j\geq2$ claiming $j-1$ or $j+1$ is
a suitable choice for Breaker, and this choice is determined by $j$
$\text{(mod 5)}$. Apart from the case $j=1$, the board always splits
into two boards of the form $F\left(a\right)$ and $G\left(b\right)$
for some $a$ and $b$ with $n=a+b+2$. We use the same notation as
in the earlier cases. \\

\textbf{Case 3.1: }$j=1$. \\

Suppose Breaker claims the element  $2$. Then $s_{i}=0$ and the
board becomes $F\left(n-2\right)$. Hence $d_{2}=\left\lfloor \frac{n+2}{5}\right\rfloor -\left\lfloor \frac{n}{5}\right\rfloor $,
which is certainly always non-negative. Since $N_{1}$ decreases by
at most $1$, it follows that $d_{1}\geq-1$. 

If $n\equiv3\text{ or }4$ then we have $d_{2}\geq1$ and hence $d_{1}+d_{2}\geq0$,
as required. Otherwise $N_{1}$ is certainly not decreasing, so $d_{1}\geq0$.
Thus we again have $d_{1}+d_{2}\geq0$, which completes the proof
of Case 3.1. \\

\textbf{Case 3.2: $j\neq1$ }and $n\equiv3\text{ or }4$. \\

The proof is identical to the proof of Case 1.1.1. \\

\textbf{Case 3.3:} $j\neq1$ and $n\not\equiv3\text{ and }4$. \\

The proof is identical to the proof of Case 1.1.2. \\

This completes the proof of Claim 3, and hence Lemma 5 holds by induction.
$\square$\\

Recall from the Introduction that $H_{b}\left(n\right)$ is the game
played on the same board as $H\left(n\right)$, but with Breaker having
the first move. Also recall that we have $u\left(P_{n}\right)=\gamma_{b}\left(n\right)$
and $u\left(C_{n}\right)=\alpha\left(n-1\right)$. We now deduce Theorems
1 and 2 from our earlier results. \\

\textbf{Proof of Theorem 1. }Note that Lemma 4 implies that $u\left(P_{n}\right)=\gamma_{b}\left(n\right)\geq\left\lfloor \frac{n+4}{5}\right\rfloor $.
In order to prove the upper bound, consider the game $H_{b}\left(n\right)$
and suppose that Breaker claims the element $n$ on her first move.
After this subsequent move, the game is equivalent to the game on
the same board as $G\left(n-1\right)$ with Maker having the first
move. Hence it follows that $\gamma_{b}\left(n\right)\leq f\left(\emptyset;\,n-1;\,\emptyset\right)$,
and thus Lemma 5 implies that $\gamma_{b}\left(n\right)\leq\left\lfloor \frac{\left(n-1\right)+5}{5}\right\rfloor $.
Therefore we have $u\left(P_{n}\right)=\left\lfloor \frac{n+4}{5}\right\rfloor $,
as required. $\square$\\

\textbf{Proof of Theorem 2. }Recall that we have $u\left(C_{n}\right)=\alpha\left(n-1\right)$.
Hence Lemma 3 implies that $u\left(C_{n}\right)\geq\left\lfloor \frac{\left(n-1\right)+2}{5}\right\rfloor $,
and Lemma 5 implies that $u\left(C_{n}\right)\leq f\left(n-1;\,\emptyset;\,\emptyset\right)=\left\lfloor \frac{\left(n-1\right)+2}{5}\right\rfloor $.
Thus it follows that $u\left(C_{n}\right)=\left\lfloor \frac{n+1}{5}\right\rfloor $,
as required. $\Square$\\

In particular, for both $G=P_{n}$ and $G=C_{n}$ it follows that
the asymptotic proportion of isolated vertices is $\frac{1}{5}$ when
both players play optimally. \\

There are many questions that are open concerning the value of $u\left(G\right)$
for general $G$. Dowden, Kang, Mikala\v{c}ki and Stojakovi\'{c} \cite{key-7}
gave bounds for $u\left(G\right)$ that depended on the degree sequence
of the graph $G$. As a consequence they concluded that if the minimum
degree of $G$ is at least $4$, then $u\left(G\right)=0$. They also
noted that there exists a $3$-regular graph with $u\left(G\right)>0$,
and they proved that the largest possible proportion of untouched
vertices among all $3$-regular is between $\frac{1}{24}$ and $\frac{1}{8}$.
It would be interesting to know what the exact value is. Their example
for the proportion $\frac{1}{24}$ is not connected, so it would also
be interesting to know what the maximal proportion is for connected
$3$-regular graphs. 

They also proved that if $T$ is a tree with $n$ vertices then we
have $\left\lceil \frac{n+2}{8}\right\rceil \leq u\left(T\right)\leq\left\lfloor \frac{n-1}{2}\right\rfloor $.
The upper bound is tight when $T$ is a star, but they did not find
a similar infinite family of examples for which the lower bound is
tight. It would be interesting to know whether this lower bound is
asymptotically correct.

\end{document}